\documentclass[letterpaper,10pt]{IEEEtran}

\usepackage{amsmath} 
\usepackage{amssymb} 
\usepackage{MnSymbol}
\usepackage{eurosym}
\usepackage{pdfsync}
\usepackage{color}
\usepackage{ulem}
\usepackage{url}
\usepackage[symbol]{footmisc}

\usepackage{url}
\usepackage{graphicx,xcolor}
\usepackage{verbatim}
\makeatletter
\newcommand{\verbatimfont}[1]{\def\verbatim@font{#1}}%
\makeatother
\verbatimfont{\ttfamily\small}

\usepackage[left,pagewise]{lineno} 
\usepackage[english]{babel} 
\usepackage{blindtext}

\title{Behavioral  Economics  for  Human-in-the-loop  Control  Systems Design\\ {\Large Overconfidence and the hot hand fallacy}}
\author{Marius Protte$^\dag$,\thanks{$^\dag$Faculty of Business Administration
    and Economics, Paderborn University, 33098 Germany. \texttt{marius.protte@uni-paderborn.de}, \texttt{rene.fahr@uni-paderborn.de}}  \and Ren\'e Fahr$^\dag$,	 \and Daniel E.\ Quevedo$^\sharp$,\thanks{$^\sharp$School of Electrical Engineering \& Robotics, Queensland University of Technology, QLD 4000, Brisbane, Australia. \texttt{dquevedo@ieee.org}} \IEEEmembership{Senior Member, ̃IEEE}}

\newif\ifPDF \ifx\pdfoutput\undefined\PDFfalse \else\ifnum\pdfoutput > 0\PDFtrue \else\PDFfalse \fi \fi
\ifPDF 
\usepackage[pdftex, plainpages = false, colorlinks=true, linkcolor=black, citecolor = green!50!blue, urlcolor = blue, filecolor=black, pagebackref=false, hypertexnames=false,  pdfpagelabels ]{hyperref}
\fi

\begin{document}
\maketitle

\begin{abstract}
Successful design of human-in-the-loop control systems requires appropriate models for human decision makers. Whilst most paradigms adopted in the control systems literature hide the (limited) decision capability of humans, in behavioral economics individual decision making and optimization processes are well-known to be affected by perceptual and behavioral biases. Our goal is to enrich control engineering with some insights from behavioral economics research through exposing such biases in control-relevant settings.

\par This paper addresses the following two key questions:
 
\begin{enumerate}
\item How do behavioral biases affect decision making?
\item What is the role played by feedback in human-in-the-loop control systems?
\end{enumerate}

Our experimental framework shows how individuals behave when faced with the task of piloting an UAV under risk and uncertainty, paralleling a real-world decision-making scenario. Our findings support the notion of humans in Cyberphysical Systems underlying behavioral biases regardless of -- or even because of -- receiving immediate outcome feedback. We observe substantial shares of drone controllers to act inefficiently through either flying excessively (overconfident) or overly conservatively (underconfident). Furthermore, we observe human-controllers to self-servingly misinterpret random sequences through being subject to a ``hot hand fallacy''. We advise control engineers to mind the human component in order not to compromise technological accomplishments through human issues.

\end{abstract}
 
\section{Introduction}
Our century is bringing interesting challenges and opportunities that derive
from the way that digital technology is shaping how we live as individuals and
as a society. A key feature of a number of engineered systems is that they
interact with humans. Rather than solely affecting humans,
often people make 
 decisions that have  an effect on the
engineered system.  As an example, when
driving our cars we often  decide to take a route which differs from that
suggested by the navigation system. This information is fed back to the
service provider and henceforth used when making route suggestions to  other
users. The analysis and design of such \textit{Cyber-physical Human Systems}
(CPHSs) would
benefit from an understanding of how humans behave.
However, given their immense complexity, it is by no means clear how to 
formulate appropriate
 models for  human decision makers, especially when operating in closed loop systems.

\par CPHSs fit into the general context of human-machine systems and pose multidisciplinary challenges that have been tackled in a number of domains. For example, interesting surveys on signal processing approaches include \cite{vemkai18} and \cite{nargeo13a}, whereas \cite{scherd13} adopts a computer science viewpoint. More focused on specific applications, the articles in the special issue \cite{tantol17} as well as \cite{dressl18a}  study the interaction between humans and vehicles.  The paper \cite{kolwal16} surveys systems comprised of humans and robot swarms and \cite{youpes20} reviews  telemanipulation by small unmanned systems.
Also within the control systems community significant advances have been made, see, e.g.,
\cite{hokspo06a,vanOverloopetal.:2015,erccar17,inogup19}. Here it is common to model the effect  of humans as a limited
actuation resource or unknown deterministic dynamics to be identified. This
opens the door to use various robust control and game theoretic
methods to design 
closed loop systems. However, most 
paradigms adopted in the systems control literature so far hide  
the decision capability of humans and limitations to their cognitive and
computational 
capabilities.  Notable exceptions to
the literature  
include \cite{lamsas14} and \cite{Fengetal:2016} where human decision making is
characterized via a Partially Observed Markov Decision Process (MDP)
\cite{puterm94} with states 
representing basic physical human aspects, such as operator fatigue or
proficiency. In particular, \cite{Fengetal:2016} investigates a stochastic two
player game, resulting from the interaction between a human operator and an
unmanned aerial vehicle (UAV). Here the UAV is allowed to react to the decisions
made by the operator, to compensate for potential non-cooperative
behaviour. 
For a CPHS with multiple human decision makers, \cite{albyil19}
models limitations in their decision making capabilities using ``level-k reasoning''
and associated game
theory \cite{camho04a}: Each agent optimizes a cost function using only reduced
knowledge of the other agents' policies.    

\par In the present work we investigate human decision making  from a 
 behavioral economics perspective. We pose the question how behavioral biases affect decision making, since we expect behavioral influences to potentially hamper technological optimization efforts.  More generally, we aim to direct attention to the human as a  poorly observed uncertainty source in CPHS. We thus consider humans 
as being 
``bounded rational'', i.e., they are strategic thinkers who want to maximize
their benefits but, at the same time, have only limited cognitive and computational
capabilities. We distinguish between closed loop scenarios, where the decision maker receives feedback about the success of its actions, and open loop situations, where no such feedback information is available. Here, we pay special attention to the role played by feedback on behavioral biases in form of
overconfidence \cite{Schaefer/Williams/Goodie/Campbell:2004} and the so-called
hot hand 
fallacy  \cite{Gilovich/Vallone/Tversky:1985}. The latter effect may arise when
a decision maker becomes aware of her past success and 
overestimates future success probabilities. 
 We will next introduce some prominent findings from
behavioral economics that can be regarded as influential for CPHSs. These will be
subsequently applied  to an experimental 
human-in-the-loop framework inspired by the UAV piloting scenario of
\cite{Fengetal:2016}. Throughout our presentation we will highlight the value of giving greater
consideration to human behavior and behavioral biases in control engineering. 
  The central outcome of our study
is that frequent feedback to human decision makers may lead to sub-optimal
results. This stands in contrast to situations wherein 
computers carry out optimizations and more information is
always beneficial \cite{bartse74,puterm94,bertse05b}.
 As we shall see, when humans act as
decision makers, then the situation is more ambiguous.       

\section{Contributions from Behavioral Economics}

Behavioral economics is a comparably young research field that is concerned with
the overarching question of how humans actually do behave in economic decisions,
in contrast to how they are prescribed to behave by economic theory
\cite{Kahneman:2003}. Therefore, traditional economic assumptions and models are
revised 
and enriched by insights from psychology. This serves to improve the understanding
and predictability of human behavior, especially regarding recurring errors and
biases in decision making \cite{Camerer/Loewenstein:2003}. 
The perception of humans thereby shifts from the assumption of a
 rationally optimizing Markov decision maker towards regarding humans as
rather bounded 
rational agents. Humans constantly use cognitive
mechanisms of simplification, namely decision heuristics, in order to process
information and make decisions under uncertainty \cite{Kahneman:2003}, \cite{Ariely:2009}. Such heuristics are used subconsciously due to individuals
not managing to adequately process the complexity of a decision problem or to
take all relevant information into account \cite{Camerer/Loewenstein:2003}. In fact, when making judgments or estimations of events, frequencies or
probabilities under uncertainty, individuals do not always obey Bayesian rules
and statistical logic, as they are meant to do in models assuming perfectly
rational agents and Markov decision makers. Such heuristic simplifications sometimes yield reasonable judgments, but may also
lead to severe and systematic errors \cite{Kahneman/Tversky:1973}, \cite{Tversky/Kahneman:1974}, \cite{Todd/Gigerenzer:2003}.  
Interestingly, even experienced researchers and professionals often underlie the
same judgment heuristics and biases as laypersons do \cite{Tversky/Kahneman:1974}.

\subsection{Overconfidence}
\label{sec:overconfidence}
Overconfidence is well-known to be one of multiple behavioral stylized facts
influencing and thereby biasing human decision making. Overconfidence can be
defined as a general miscalibration in beliefs \cite{Lichtenstein/Fischhoff:1977}, more specifically, the discrepancy between confidence and accuracy
\cite{Schaefer/Williams/Goodie/Campbell:2004}. In being overconfident, people overestimate their
own capabilities, their knowledge or the general favorability of future
prospects \cite{Barber/Odean:2001}. The authors of \cite{DeBondt/Thaler:1995} go as far as labeling
overconfidence the “most robust finding in the psychology of judgment”. 
 
Early contributions to research on overconfidence found that people
systematically tended to be unrealistically optimistic about their future, in
judging themselves to be more likely to experience a variety of positive events,
and to be less likely to experience negative events, compared to others. This
pattern has been traced back mainly to the degree of desirability influencing
the perceived probability of such events \cite{Weinstein:1980}. In the original
study  \cite{Svenson:1981} over $80\%$ of the subjects regarded themselves as
more skillful and less risky car drivers than the average driver. Moreover, one
half of the subjects estimated themselves to be amongst the top $20\%$ of the
sample, vividly illustrating overconfidence as the discrepancy between
confidence and accuracy (belief and reality) in doing so. 
Meanwhile, evidence on overconfidence does not depend on whether or not
experimental subjects are familiar to the given tasks or actions \cite{Hoelzl/Rustichini:2005}. The task’s difficulty, on the other hand, is known to
have an impact on overconfidence, as more difficult tasks were shown to
facilitate overconfidence rather than easy tasks \cite{Lichtenstein/Fischhoff:1977}, \cite{Moore/Healy:2008}.
Overconfidence can express itself, and consequently has been studied, in multiple
ways, like `better-than-average' beliefs, e.g. \cite{Hoelzl/Rustichini:2005}), or `overprecision', e.g. \cite{Radzevick/Moore:2011}). The facet of overconfidence most fitting for the
purpose of our current work appears to be overestimation, described as self-servingly
overestimating the likelihood of desirable outcomes, supposedly fueled by
wishful thinking for such desirable outcomes to occur\footnote{Overconfidence is
  therefore regarded as closely related to self-deception, see e.g. \cite{vanHippel/Trivers:2011}. 
  }  \cite{Moore/Schatz:2017}, \cite{Mayraz:2011}.  
An experimental approach for
studying overconfidence was requested in \cite{Hoelzl/Rustichini:2005}, as most
research in this field merely relies on verbal 
statements or subjective estimations, rather than monitoring human
decision-behavior. Our current analysis makes a contribution to this
call. Overconfidence is connected to a variety of other behavioral phenomena,
which will be revised in the sequel.

\subsection{Underestimation of systematic risks}
\label{sec:systematic}
Analogous to overestimating favorable outcomes, overconfidence also implies underestimating risks, or rather the variance of risky processes,
indicating a too narrow distribution in one's subjective probability beliefs \cite{BenDavid/Graham/Harvey:2013}. Often, humans  
assess probabilities incorrectly or rather draw incorrect conclusions from them through not adhering to Bayesian rules or neglecting base-rates \cite{Kahneman/Tversky:1973}, \cite{Sedlmeier/Gigerenzer:2001}. In doing so, smaller probabilities are usually overestimated while larger probabilities are underestimated \cite{Kahneman/Tversky:1979}, \cite{Slovic/Fischhoff/Lichtenstein:1982}. Vice versa, decision makers are also found to underestimate
the overall likelihood of the occurrence of risks with small single
probabilities of occurrence in some instances. Especially risks of disqualification
\cite{Abbink/Irlenbusch/Renner:2002} through low-probability events generate less concern than their probability warrants on
average \cite{Weber:2006}. A cognitive
process behind this underestimation of especially very low probability risks, like
being involved in car crashes or natural disasters, can be characterized as
subconsciously approximating low probabilities with zero in order to not having
to mentally deal with them anymore \cite{Slovic/Fischhoff/Lichtenstein:1978}.  

\par Consequences of this issue can be observed in various areas of economic and
social decision making. Overconfident individuals, who underestimate health,
financial or driving risks, are more likely to have insufficiently low insurance
coverage against such risks \cite{Sandroni/Squintani:2013}. Furthermore,
employees are observed to underestimate the risk of their own company's stock
and being overly optimistic about its future performance. They   tend to
include such stocks too heavily into their retirement plans, despite the
respective stock being riskier than the overall market, as a consequence of excessive
extrapolation of positive past performance \cite{Benartzi:2001}. 
Money being at stake does neither change overconfident behavior, nor does it
lead to better estimation results concerning abilities, probabilities or
risk. Experimental evidence was obtained on overconfidence leading to overly
rushed market entries, that often precede business failures, due to
entrepreneurs overestimating their relative chances to succeed \cite{Camerer/Lovallo:1999}. Overconfident top-managers were further shown to underestimate
the volatility of cash flows of S\&P 500 companies, resulting in erroneous
investment and financing decisions, in a large-scale survey over 10 years \cite{BenDavid/Graham/Harvey:2013}. Especially managers competing for
leadership positions display overconfidence by tending to take on riskier
projects due to underestimating their risks \cite{Goel/Thakor:2008}. Overconfident CEOs also overestimate the positive impact of their
leadership and their ability to successfully complete a merger to generate
future company value. These CEOs are found to execute value-destroying mergers,
thereby exposing the affected companies to high financial or even existential
risks \cite{Malmendier/Tate:2005}. Multiple experimental studies found that
overconfident investors trade more excessively than others, overestimate the
impact of little information they have \cite{Odean:1999}, \cite{Barber/Odean:2000}, and do not adjust their trading volume to new negative information
\cite{Trinugroho/Sembel:2011}. These findings, which are inconsistent with the
Markov assumption and Bayesian updating, contradict theoretical expectations of
rationally optimizing economic agents  and re-emphasize the
notion of individuals as bounded rational decision makers \cite{Kahneman:2003}.

\subsection{Attribution theory}
\label{sec:attribution}
Another mechanism contributing overconfident decision making is a biased
perception of causality in a self-serving way. Accordingly, individuals express
the tendency to attribute positive past outcomes to themselves and their
abilities, while blaming negative outcomes on external circumstances
\cite{Frieze/Weiner:1971}. Successes are internalized, failures are externalized
\cite{Langer:1975}. In an early experiment, teachers attributed learning
successes of their pupils to their teaching skills, while blaming bad learning
performances on the pupils themselves \cite{Johnson/Feigenbaum/Weiby:1964}. 
 This self-serving attribution bias was identified to be the major driver of
the aforementioned CEO overconfidence leading to risky mergers \cite{Billett/Qian:2008}. In \cite{Gervais/Odean:2001} overconfidence evolves when traders get feedback about their ability through experience in a multiperiod market model and overweigh the role of their ability on prior success. Overconfident investors subsequently tend to become even
more overconfident in their future investments. Meanwhile, financial losses are
rather attributed to environmental circumstances, like unfavorable macroeconomic
developments or simply bad luck, with the investor’s degree of overconfidence
remaining constant \cite{Hilary/Menzly:2006},  \cite{Daniel/Hirshleifer:2015}.
Also, professionals were found to be as likely as laypersons to express
overconfidence in making economic decisions as well as in re-evaluating the
quality of their own previous decision in hindsight \cite{Torngren/Montgomery:2004}. 

People overestimate their own capabilities as well as their control over future
events that are actually determined by external factors, especially chance. This
was illustrated in an experimental study in which participants were either given
a lottery ticket (control group) or allowed to select a lottery ticket
themselves (treatment group). Subjects who had chosen their ticket themselves
were, on average, demanding four times more money than those who were simply
given their tickets, when asked to name a price they would sell it for. Those
who chose the ticket themselves mistakenly assumed they would therefore be more
likely to win, although the winners would be determined entirely by chance. This
mechanism has been labeled the ``illusion of control'' \cite{Langer:1975},
although critics of this concept have argued that it rather identifies a pattern of
people overestimating their ability to predict future outcomes than of
people overestimating their ability to control future outcomes \cite{Presson/Benassi:1996}. Both
interpretations are regarded as sufficient for the purpose of our current analysis.

\subsection{The role of feedback in overconfidence}
\label{sec:feedback}
Providing feedback is commonly regarded as an efficient learning mechanism in
Markov Decision Processes in order to compute optimal behavior
\cite{vanOtterlo/Wiering:2012,puterm94,bartse74} and offers great value for deriving
logical inference rules in accordance with Bayesian reasoning in machine
learning  \cite{Kasnecietal:2010,bertse19a}. While these points may hold true
for cyber-physical systems without human decision makers, humans are often
observed to react to feedback 
differently. Contradicting Markovian assumptions, behavioral economic research
constantly observes errors in updating the information set (e.g.,
\cite{Camerer:1999}, \cite{Charness/Levin:2005}). This indicates that providing
more feedback information to humans does not necessarily lead to better
estimations and decisions. The importance of outcome feedback in counteracting
certain behavioral biases, through continuous ``monitoring of progress through
a judgment-action-outcome loop''  \cite{Kleinmuntz:1985} and offering corrective
adjustments, has been pointed out before  \cite{Hogarth:1981}, but it may not
apply equally to every type of bias. 

\par While research studying the connection of high information supply to
overconfidence is comparably scarce, in various instances a higher amount of feedback was not found to be able to improve
behavioral rationality with regards to overconfidence, as it would be expected
by theory. Both overconfidence and underconfidence were shown to persist in employees' choices of incentive schemes although receiving clear feedback revealing the optimal alternative for them \cite{Larkin/Leider:2012}. Experiments on sports outcome-estimation
showed that overconfidence did not decrease but rather increase with additional
information, as the subjects’ confidence increases more than their accuracy does
\cite{Tsai/Klayman/Hastie:2008}. Similarly, CEO overconfidence in forecasts was found
to persist against corrective feedback \cite{Chen/Crossland/Luo:2015} and
venture capitalists were shown to be more overconfident when having access to
more information \cite{Zacharakis/Shepherd:2001}. In general knowledge tasks that
are related to each other, corrective outcome feedback has shown able to
mitigate underconfidence but not overconfidence \cite{Subbotin:1996}. It turns
out that feedback
appears to affect behavioral biases rather ambiguously with overconfidence
tending to largely persist, or even growing, in response to additional
information.

\subsection{Misperceptions of random sequences}
With regards to human sequential decision processes, another
potential threat of high 
information supply leading to  biased decision making, is the danger of falling subject to
the “hot hand fallacy”. The study of this error in the assessment of past
outcomes originates in observations from basketball. If this player
repeatedly scored on his past throws, he is commonly judged to have a “hot hand”
by his teammates and the audience. Therefore, the player's likelihood of scoring
on future throws is judged to increase by the audience, even though future throws
are independent of past throws. In doing so, the player's objective
probability of scoring is overestimated based on his past successes
\cite{Gilovich/Vallone/Tversky:1985}. This hot hand fallacy results from
misjudging statistically independent favorable events to be connected to one
another, implying a positive autocorrelation between them. Repeated past
positive outcomes are therefore erroneously expected to occur again in the near
future, thereby overestimating their objective probability of occurrence. The
hot hand fallacy is a very vivid illustration of biased probability judgments
through presumed representativeness of recent information. Misinterpretations of
random sequences through such extrapolation of recent outcomes into the future
were empirically found to apply to several contexts. Regarding basketball, the
findings by \cite{Gilovich/Vallone/Tversky:1985} were supported by studies on
basketball betting market odds, although only small effects were found
\cite{Camerer:1989}, \cite{Brown/Sauer:1993}. In a laboratory experiment,
simulating a blackjack game, gamblers were observed to bet more money after a
series of wins than after a series of losses \cite{Chau/Phillips:1995}. In
more sensitive contexts, individual decisions show signs of the hot hand fallacy
as well. For instance, financial investors who had negative prior experiences
with low stock returns, exhibit a decreased willingness to take financial risks
in future investments \cite{Malmendier/Nagel:2011}. 
Attribution theory, especially the illusion of control mentioned above, may
factor into the hot 
hand assumption as well, since subjects are more likely to attribute recent
sequences with low alternation to human skilled performance, while sequences
with high alternation are rather perceived as chance processes  \cite{Ayton/Fischer:2004}. 

Before proceeding, we note that the hot hand fallacy's opposing bias, the so-called \textit{gambler’s fallacy}, presents the erroneous expectation of systematic reversals in stochastic processes. These are grounded in the belief that a small sample should be representative for the underlying population, i.e., that a Bernoulli random process should balance out even across few rounds \cite{Rabin/Vayanos:2010}.
Both phenomena have recently been strongly disputed in their status as biases \cite{Miller/Sanjurjo:2018}, \cite{Miller/Sanjurjo:2018b}.

\begin{figure}[tp]
\centering
\includegraphics[scale=.25]{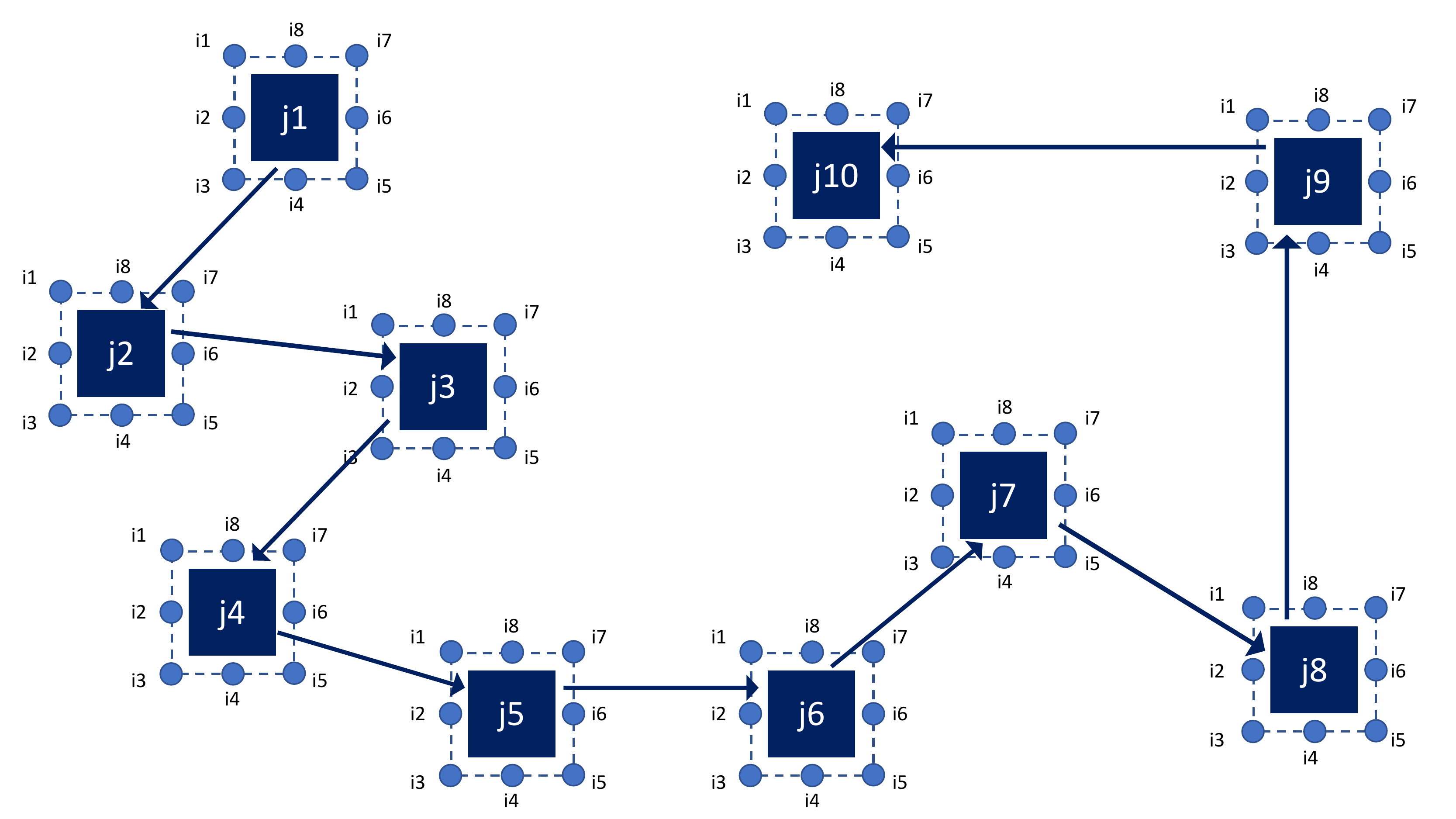} 
\caption{Illustration of a road network for UAV missions, adapted from
  \cite{Fengetal:2016}. A single drone follows a given path overflying $10$ junctions. At each junction $j\in\{1,2,\dots 10\}$ a maximum of $N=8$ rounds can be flown.}
\label{fig:dronepath}
\end{figure}

\section{Surveillance drone piloting framework and hypotheses}
\label{sec:setup}
Our experimental design for analyzing behavioral biases in CPHSs builds upon a
basic setup involving UAVs 
that interact with human operators, as discussed in \cite{Fengetal:2016} and which
is graphically illustrated in Fig.\ \ref{fig:dronepath}.  Partly autonomous UAVs
are used for 
road network surveillance 
inside a network that connects multiple traffic junctions. The drone pilot's
objective is to gain the maximum information about the traffic at all  
junctions using a single drone which follows a predesigned path, see Fig.~\ref{fig:dronepath}. To increase the picture quality (and thereby the information obtained), the
operator may choose to fly up to 8 rounds over each junction $j$, each round taking an additional photo. The individual photos are combined into a higher-resolution 
image, whose total value will be denoted $I_j$. 
This total value depends on the number of flights, say $f_j$, and on the random additional information content of each photo taken. The situation is modelled via the random process $\sigma_j(i)$ which quantifies the information content of the combined photo of junction $j$ at round $i$. The process is initiated at $\sigma_j(1) =25$ and through iteration of
\begin{equation}
\label{eq:foto}
    \sigma_j(i+1)=\sigma_j(i) + s_j(i)\rho(\sigma_j(i)), \quad i=1,2,\dots,f_j-1,
\end{equation}
 leads to $I_j = \sigma_j(f_j)$. 
\par Here, the  additional information content of the photo taken at round $i$ is quantified by the concave function $\rho(\cdot)$ which describes a decreasing marginal yield:
\begin{multline}
\label{eq:rho}
    \rho(25)=25,\; \rho(50)=20,\; \rho(70)=10, \\ \rho(80)=\rho(85)=\rho(90)=\rho(95)=5.
\end{multline}
The above reflects the fact that each additional picture often comes with less information value added, since a rough image of the traffic flow is already gained by the earlier pictures, and subsequent pictures only help further sharpening the image.
 In \eqref{eq:foto}, $s_j(i)\in\{0,1\}$ is an i.i.d.\ Bernoulli random process, with probability $P[s_j(i) =1] =p$. This models instances wherein the new picture offers no information value added over the previous one due to, for instance, bad weather causing poor visibility, strong wind or  lacking flow of traffic. 
 

\par  Whilst taking more photos will often lead to combined images of higher quality and also with more information content, after  each individual flight there exists a (small) probability of $r$
of the drone crashing,  leading to a loss of $D$ (the value of the UAV) and the inability to continue flying again over the current or later junctions. 

\par Given the above, the  value of the  image obtained at each junction $j$ belongs to the finite set
$\{0,25,50,70,80,85,90,95,100\}$.
The current combined value of the drone and images after flying   $i$ rounds at junction $j$ satisfies
\begin{equation}
\label{eq:Vji}
    V_j(i) = D c_j(i) + \sigma_j(i) + \sum_{\ell=1}^{j-1}I_\ell, 
\end{equation}
where
\begin{equation}
\label{eq:cji}
    c_j(i)  = \begin{cases}
    1, &\textrm{if the drone is still intact}\\
     &\textrm{$\quad$after flying   $i$ rounds at junction $j$,}\\
0, &\textrm{otherwise.}
\end{cases}
\end{equation}
The  total value gained by the operator at the end of the mission is  given by:
\begin{equation}
\label{eq:value}
V=400 c_{10}(f_{10}) + \sum_{j=1}^{10}I_j.
\end{equation}

\begin{figure}[tp]
\centering
\includegraphics[scale=.25]{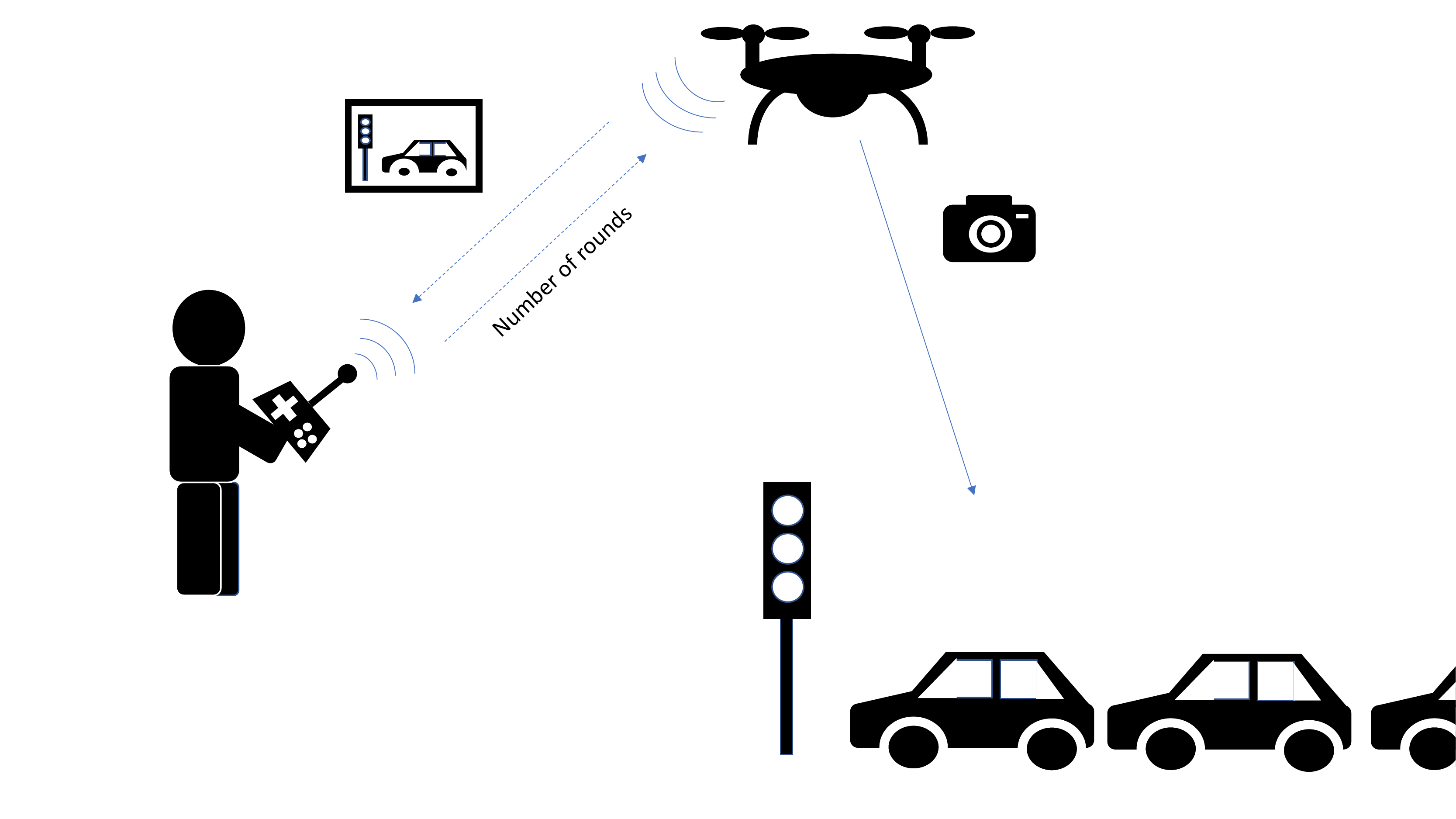} 
\caption{Human-in-the loop control system: The operator decides how many times the UAV should fly over each junction  in order to gain accurate information about traffic conditions in the city. To assist in the decision making, the drone may send feedback information about the picture quality, i.e., information gained.}
\label{fig:trafficjunction}
\end{figure}

\subsection{Decision problem}
Given the possibility of the drone crashing, rather than simply intending to fly the maximum number of 8 rounds over each junction, the operator is faced with the decision of 
how many rounds to fly over each junction in order to maximize its profit $V$. This amounts to a sequential stopping rule problem defined over a finite horizon.
 To make decisions when to fly (or switch) to the next junction, in addition to knowing the system model and its probabilities, the operator receives feedback from the drone, see Fig.~\ref{fig:trafficjunction}.  Every time the UAV passes over a junction, it
sends back $\sigma_j(i)$, the  information gain obtained so far over the current junction. Further the operator, is informed whether the drone has crashed. Using this  information, the
human operator is faced with the task of designing a closed loop flight policy.  
\par
\begin{figure*}[tp]
\centering
\includegraphics[scale=0.55]{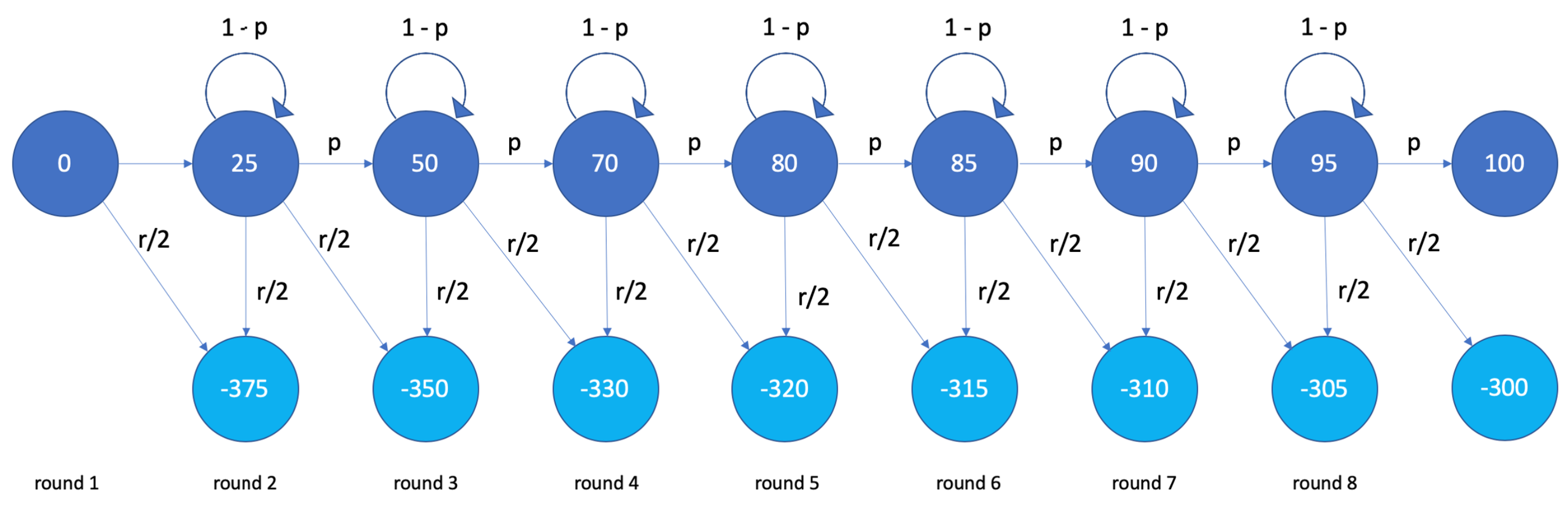} 
\caption{Marginal gains within each junction, for $D=400$. Diagonal arrows describe an increase in information value with a
  subsequent crash, vertical arrows represent a non-increase with a subsequent
  crash. Note that the maximum  information value
gain attainable at each traffic junction, namely 100, can only be reached if a
subject decides to fly all eight 
 rounds possible while obtaining an increase in each round.  
} 
\label{fig:qualitygains}
\end{figure*}

 The exact stopping rule can, in principle, be derived using dynamic programming~\cite{bertse05b,fergus08}.  
However, the value function depends not only on the  value $\sigma_j(i)$, but also on the current round $i$ as well as the junction $j$, with later junctions having less value. Thus, instead of pursuing an optimal strategy, a 1-stage look-ahead rule within the current junction becomes a reasonable alternative. Using such a myopic policy, the operator chooses to flight another round over the current junction $j$ if and only if the drone is intact ($c_j(i)=1$), $i\leq 8$ and  the marginal gain of flying one more round is positive:
\begin{equation}
   \label{eq:drift}
   g(\sigma) \triangleq E[V_j(i+1)-V_j(i)\,|\, \sigma_j(i)=\sigma] =p\rho(\sigma)- D r>0,
\end{equation}
where we have used \eqref{eq:Vji}, see also Fig.~\ref{fig:qualitygains}.

\par In the sequel we will fix the parameters to $D=400$, $p=0.5$ and $r=0.02$, so that \eqref{eq:rho} provides
\begin{multline}
g(25)=4.5,\; g(50)=2, \; g(70)=-3, \\ g(80)=g(85)=g(90)=g(95)=-5.5.     
\end{multline}

This  makes
flying  until an information value $\sigma_j(i)=70$ is accumulated in the current junction
the optimal choice. From the fourth node onwards, the marginal gain from flying
another round is lower than the marginal loss, therefore a rational (and risk-neutral) agent should
refrain from flying further rounds.

\par In addition
to this full feedback case, we will also consider a situation wherein the drone
does not transmit the feedback information $\sigma_j(i)$ to the operator. This requires the
operator to design an 
open loop policy \cite{bartse74}, as detailed below.

\subsection{Experimental drone framework}
\label{sec:design}
In our behavioral economic engineering project (see Sidebar ``Behavioral Economic Engineering"), we abstract
the framework into an economic experiment (see Sidebar ``Experimental research method and Induced Value theory") in order to gain
insights on actual individual behavior in human-in-the-loop control. The experiment translates the CPHS concepts of  open versus closed loop control into a sequential decision making game, framed in a context
in which the subjects are owners to an UAV with a photo function. They
are told that they were hired by the local city administration to support the
city's traffic surveillance by taking photos of ten important traffic
junctions. As outlined above, the subjects' task is to decide, how many rounds the UAV should fly
over each of the 10 traffic junctions (up to a maximum of 8 times), while the drone would fly
and take pictures autonomously. Subjects are incentivized, as they are paid
according to the amount of information value gained through the pictures their
drones take in the fictitious currency ``Taler'', with one Taler equaling one
unit of information value gathered. Taler are exchanged into Euro, after all
parts of the experiment are completed. However, during each round in which a UAV
 flies, it faces a constant crash risk of $r=2\%$. At the end of this job,
i.e., after all ten traffic junctions, subjects are able sell their drone for an additional
fixed payment of $D=400$ Taler, as long as it is still intact at this point, see   \eqref{eq:value}. 

\par Each traffic junction represents one stage of the experiment.  For
the sake of 
simplicity, subjects do not have to bear any costs from flying the
drone. Further, the drone's value does not diminish from usage and subjects are
told that batteries for the drones as well as (ground) transport between the traffic
junctions is provided free of charge by the city administration. One battery charging allows
the drone to fly up to eight rounds before the drone must land to recharge. Therefore, subjects are able
to decide for a maximum of eight rounds per traffic junction and
consequently take up to eight pictures. This corresponds to a maximum of eight
 nodes that can be reached per junction. In the first round flown over each traffic junction, an information value of $\sigma_j(1)=25$ is gained with certainty, as any picture
taken represents an increase in information value over no picture taken at
all. From the second round onwards, the information value potentially gained per round
decreases,  see \eqref{eq:rho}.  For the sake of
 simplification, a crash can only happen after a picture is taken. As per~\eqref{eq:value}, the
 information value gained at the junction where the crash occurs up to the point
 of it, still qualifies for the payment, as subjects were  told that the
 drone's  memory chip would survive the crash.

\par The operators need to  decide how many rounds a UAV should fly over
each traffic junction in order to maximize the total gain $V$ in~\eqref{eq:value} and 
consequent payoff. 
In order to do so, some of them receive feedback on
the results of the previous rounds -- depending on the treatment -- which they can
base their decision on. 
As noted before, the resulting optimization problem belongs to the class of general stopping problems. Its solution, both with and without feedback information, can in principle be
derived, but requires careful computations.  Since individuals are well-known to have limited
cognitive and computational capabilities, they will not be able to compute a globally
optimal strategy.   
Instead, optimizing subjects will use the decision heuristic in~\eqref{eq:drift}, i.e., the fourth node is identified as the one to
reach.

\par The above leads to two decision heuristics, depending on whether
feedback information is available or not. In the former (\textbf{closed loop heuristic}),
at every junction 
the UAV should aim to fly as many rounds as 
needed until reaching an information value of $\sigma_j(i)=70$.  To be more precise,   introduce 
$a_j(i)\in\{0,1\}$, where $a_j(i)=1$ corresponds to the decision to fly another round over the current junction, then the closed loop heuristic	 can be stated via: 
\begin{multline}
\label{eq:cl}
a_j(i) = 1,\quad  \textrm{if and only if  $\sigma_j(i) < 70$ and $i\leq 8$} \\ \textrm{and the drone is intact.}
\end{multline}
Note that in this case, the total number of rounds flown over a junction $j$, namely $f_j=\sum a_j(i)$, depends on
the sequence of increases and non-increases the subject experiences. 
\par  In the
absence of feedback from the drone, a suitable  \textbf{open loop heuristic}
amounts to  attempting to fly $f_j=5$ times at every junction:
\begin{equation}
\label{eq:ol}
a_j(i) = 1,\quad  \textrm{if and only if $i\leq 4$ and the drone is intact.}
\end{equation}
 This  yields the
expected result of three increases in information value, with the first one
being granted with certainty in the first round.
The above rules of thumb present reasonable
heuristic approximations to the task's optimal solution when implemented for all
junctions.
A rational and risk-neutral decision maker, who acts
according to expected utility theory \cite{vonNeumann/Morgenstern:1944}, would follow
these heuristics. Deviations from these decision-rules upwards or downwards represent
indicators of overconfidence or underconfidence respectively, depending on the
direction.

\par The experiment ends either after all junctions are completed and the drone is sold
or once the drone has crashed. A crash could be caused by external factors like
weather conditions, the drone hitting some object or animal, or technical
failure. This risk of a total loss of the drone's value in case of a crash
imposes the prospect of high one-time costs onto the participants and thus
prevents them from simply choosing the maximum number of rounds possible at every
traffic junction. As subjects are paid according to the information value gained, striving
for the optimal information value (using the above heuristic) presents the
rational strategy. Always flying 
the maximum number of 8 rounds is not efficient, as gaining information value has
to be weighed against maintaining the chance of being able to sell the drone at the
end of the experiment, as indicated in \eqref{eq:value}.

\subsection{Hypotheses}
\label{sec:hyp}
The variety of behavioral phenomena and stylized facts regarding decision biases
and errors in relation to overconfidence presented in the preceding section,
suggests inefficient human decision making. In fact, individuals are commonly observed to
act overconfidently. In doing so, they tend to underestimate systematic risks,
insufficiently process joint and conditional probabilities, self-servingly
misattribute causal relationships and erroneously expect favorable outcome
sequences to continue in the future. These aspects should be particularly
relevant in our example, since the probability of the drone crashing may appear
small for one specific round ($2\%$), while the overall risk of the UAV
crashing at any time   over the course of the experiment is substantially
higher. Modeling the crash probability as a Bernoulli process, the chance of a
crash when  attempting to fly the maximum number of rounds at each of the ten
junctions would be   $1-0.98^{80} > 0.8$. A narrow focus on the low single-decision crash 
probability \cite{Read/Loewenstein/Rabin:1999} could therefore lead subjects to
underestimate the overall crash probability and consequently take higher
risk. Providing feedback is commonly assumed to diminish these behavioral biases
but known to not being able to entirely resolve them, and in some instances even
aggravating them. Subjects may rather attribute positive feedback to their own
decision performance to themselves instead of realizing that it is caused by
chance. All these aspects point in the direction of expecting overconfident
drone piloting. 
Manager overconfidence in economic decision making can be summarized to
originate in either overestimating expected cash flows or underestimating risk
\cite{Hirshleifer/Low/Teoh:2012}. Translated to our design, the drone pilots'
overconfidence could either result from overestimating the likelihood of
information value improvements or from underestimating the risk of the UAV
crashing. While the task is framed in a quite understandable context, it is
non-trivial in its solution. Most subjects are expected to approach the task
rather intuitively, i.e., not solving it systematically through comparing
marginal gain and marginal costs, and consequently act overconfidently
through one of its facets described above, if not through both.  
According to the general consensus of literature on overconfidence in human
decision making, the following hypothesis is posed:  

\paragraph{\bf{H1}} \textit{Drone pilots will generally act overconfidently in
  flying more than the optimal number of rounds, regardless of feedback.} 
  
\par As is also known from the literature, individuals are prone to unreflectively repeat decisions
that   previously yielded them positive outcomes, even though such past successes
entirely depended on chance and have no impact on future outcomes at all. This
decision heuristic,  the hot hand fallacy already discussed, presents a major
misinterpretation of random sequences. While the hot hand fallacy in regard to
the original study on basketball throws \cite{Gilovich/Vallone/Tversky:1985} is
subject to criticism today (with the data being re-analyzed and the conclusion being reversed \cite{Miller/Sanjurjo:2018}), the general idea of individuals misjudging
the meaning of random outcome sequences remains to be tested in the context of CPHS. For sequential
decision problems featuring feedback from the system, decision makers appear especially vulnerable to fall
victim to the hot hand fallacy. Translated to our drone framework, it is
therefore expected that those drone pilots who immediately reach the optimal information value of a junction through experiencing gains in repeated rounds in the closed loop system (as the possibility
is only given there), perceive themselves to be ``on a roll'', i.e., as having an
equivalent to a hot hand\footnote{In our experiment, we define a hot hand as three consequent increases from the very start of a junction. With the junctions being distinct from each other and the first success being certain anyway, there is no possibility of pattern overlays \cite{Miller/Sanjurjo:2018b}. Therefore, as well as due to the outcome of a fourth flight being unimportant for judging the decision inefficient, we do not consider our results subject to "streak selection bias" \cite{Miller/Sanjurjo:2018}.}. Since decision makers are tempted to expect such
apparent trends to continue, they will likely fly beyond the optimal number of
rounds, overconfidently making their operation inefficient with regards to risk
and reward. This way, subjects fall victim to the hot hand fallacy\footnote{Literature observes a shift from hot hand beliefs for shorter outcome streaks towards a gambler’s fallacy pattern for increasing streak length \cite{Jarvik:1951},   \cite{Rabin/Vayanos:2010}. Drone controllers who experience multiple increases in a row would therefore begin to overestimate the probability of a non-increase at some point and stop flying to prevent a crash accordingly. }
Compared to the gambler’s fallacy, a hot hand is considered more suitable to portray drone controller’s behavior in his attempt to gather information value, with humans more prone to expect favorable outcomes to repeat under their apparent control \cite{Ayton/Fischer:2004}.\footnote{Also, a hot hand fallacy can be considered the more problematic error compared to the gambler’s fallacy, since in the case of the latter, subjects would stop earlier and not risk a drone crash through flying in periods that yield marginal losses (although sacrificing marginal gain). }. This leads to our second hypothesis: 

\paragraph{\bf{H2}}\textit{Sequences of immediate positive feedback will lead to drone pilots falling victim to the hot hand fallacy.}

These two hypotheses were tested in an experiment involving students at
Paderborn University in September 2019. Details are given next.

\section{Experimental design}
\label{sec:expdesign}
In order to be able to compare behavioral effects in closed and open loop
operation, a between-subject design experiment is conducted, meaning each
subject can only participate in exactly one treatment with the groups being
compared afterwards. 

\subsection{Treatments}
\label{sec:treatments}
 To evaluate the effect of feedback in human-in-the-loop control, two treatments
 are presented:  {Closed Loop} and  {Open Loop}. The groups differ
 only in the fact that subjects in one of them receive immediate outcome
 feedback after every decision made, while subjects in the other group do not 
 receive any feedback until the end of the experiment. 

\par In the \textbf{Closed Loop Treatment}, subjects receive feedback directly
after each round of each junction on whether an increase in information value was
gained in this specific round, see Fig. \ref{fig:qualitygains}. They were also
informed about the current total information value gained for the respective
junction and whether the drone was still intact. Subsequently, they were asked
whether they wanted to fly another round. The framing in the instructions said
that the drone would transmit the pictures taken to the drone pilot’s laptops immediately, so feedback was given just in time for the next round’s decision.
 
\par In the \textbf{Open Loop Treatment}, subjects are provided with
  no feedback at all. Subjects make decisions on how many rounds the drone should fly over each junction entirely in advance. Afterwards, they are informed about the total information value they obtained and whether
the drone is still intact and sold at the end of the experiment. The framing in
the instructions states that 
the UAV has to be deconstructed in a complex procedure to be able to extract
and read out the memory chip to check on the pictures taken. Therefore, no feedback on
the pictures  would be provided until either all junctions were completed or
the drone had crashed.

\subsection{Procedure}
\label{sec:procedure}

The experiment was computerized using the experiment
software oTree \cite{Chen/Schonger/Wickens:2016} and hosted centrally on a university
server, so subjects could remotely access the experiment and did not have to
come to the laboratory physically. 500 subjects who previously had enrolled
voluntarily into the BaER-Lab\footnote{BaER-Lab stands for ``Business and
  Economic Research Laboratory''. For further information visit
  www.baer-lab.org.}  student participant pool were chosen randomly, 
with 250 each being randomly assigned to the two treatments in
advance. These subjects were contacted via the online recruitment system ORSEE
\cite{Greiner:2015} and invited to participate within the following five days.
\par The invitation email included a hyperlink that directed the subjects to
their respective treatment, where they received the detailed
instructions for the experiment. The instructions of both treatments
only differed regarding the provision of feedback, as explained above. To be able
to progress to the drone flying task, subjects had to correctly solve four
control questions. The questions revolved around the
central parameters of the experimental design in order to assure that the
subjects had read and understood the instructions in its critical
parts. Subsequently, subjects advanced to the drone flying task, in which they
had to make their decisions through clicking single choice buttons. Subjects would receive one unit of the fictitious currency Taler
for each unit of information value gained as pay for flying the drone. During
the whole experiment, one subject's payoffs did not depend on the decisions of
any other subject. After completing the drone flying task, a result screen
presented the total information value accrued over all traffic junctions, the state of
the drone, i.e., whether it was still intact and consequently sold to earn 400
additional Taler or not, as well as the corresponding payoffs the subject generated in
Euros. Payoffs were translated from Taler to Euro at an exchange rate of 1 Euro
per 120 Taler.  
Subjects were then asked to answer a standardized questionnaire, that included questions on demographics,
perceived task difficulty and the subject's reasons for choosing the numbers of
rounds the way they did. Afterwards, subjects were presented with the multiple
price list of \cite{Dohmen/Falk/Huffman/Sunde:2010} (see Sidebar ``Multiple
Price List Format''), in order to measure their risk preferences.

\par Filling out the price list was incentivized, as for every fifteenth subject one
out of the list's twenty rows was randomly selected and the payoffs that the
subject’s chosen alternative row yielded in this row was added to the his
total payoffs. In case Alternative A had been selected by the subject, the fixed
amount stated in the respective row was simply added to its payoffs. In case
Alternative B had been selected, the lottery was automatically played out and
either \euro{0} or \euro{30} were added to the subject's payoffs. After completing the
price list, subjects were informed about whether the list came into effect for
their payoffs, and in case it did which row had been selected. On a final screen,
subjects were informed about their total payoffs in Euros and
thanked for their participation. 

\par A subject's total payoff function can consequently be formalized as
\begin{equation}
\label{eq:payoff} 
\Delta m = \frac{1}{120} V + d(\ell)_{i=15x}
\end{equation}
where $V$ is the total value, as in \eqref{eq:value}. The term $d(\ell)$ represents the
additional payoff from the multiple price list, conditional on the selected row
$\ell$ of the list that a subject was paid in case of being a fifteenth
participant, i.e., a subjects participant's ID being a multiple of 15. 

\par Subjects were able to collect their payoffs in cash at the Chair's secretariat following their participation by stating a unique
eight-figure ID-Code they had to create in the beginning of the
experiment. That way the correct payoffs could be handed out to each subject,
while maintaining anonymity about the subject's decisions made.

\section{Experimental Results}
\label{sec:experimentalresults}
Out of the 105 subjects who finished the online experiment, 57 participated in
the Closed Loop and 48 in the Open Loop treatment
respectively. The average age of the subjects was 23.87 years and varied by about
half a year between treatments. 31 subjects (29.52$\%$) were female, while one
subject classified itself as ``diverse''. On average, subjects earned \euro{4.89}
from the experiment, plus the amount they would earn from filling out the
multiple price list in case they were paid due to being a fifteenth
participant. 
 
\par Subjects in the Closed Loop treatment flew, on average, 5.89 rounds over
each traffic junction,  thereby exceeding the average number of rounds from the Open Loop treatment by 0.79 rounds per junction. Over all
rounds flown, subjects in the Closed Loop treatment gained a cumulative
information value of 651.84, on average. This represents a $27.83\%$ surplus compared to the average information value gained in the Open Loop treatment.
Similarly, a higher drone
crash rate of nearly $70\%$ was observed in the Open Loop group. These two
dimensions generally interact, since the information value is fixed, once a drone
crashes. Averages of rounds flown and accumulated information value 
as well as average earnings and crash rates for each treatment are displayed in Table \ref{tab:distrib}.

\begin{table}[tp]
\caption{Average values of rounds flown, gained information value $V$,  
  earnings $\Delta m$ and crash rates, 
  by treatment.  Standard deviations are presented in parenthesis.} 
\label{tab:distrib}
\begin{tabular}{| c | c | c | }
\hline
\hline
   &  Closed Loop  & Open Loop \\
\hline 
Average rounds flown &   5.89  &   5.10 \\
 & (2.08) & (1.82)\\
\hline
Average information value &  651.84  & 509.90 \\
gained in total & (397.56) & (380.13)\\
\hline
Average earnings & 5.43 & 4.30  \\
in Euro & (3.31) & (3.17)  \\
\hline
Crash rate & 52.63$\%$ & 68.75$\%$  \\
\hline
\hline
\end{tabular}
\end{table}

As could be expected from control theoretic results \cite{bartse74}, the Closed
Loop system is clearly 
observed to be more effective than the Open Loop piloting system, with regards
to accumulating information about traffic conditions in our
framework. Considerably more information value was aggregated with a closed loop
of immediate outcome feedback being implemented, which also translates to a
higher monetary payoff for the subjects. 
%
However, subjects displayed behavioral biases as we will see. Indeed, the average number of rounds in the Closed Loop treatment being nearly
six already hints that several subjects tended to fly beyond the
heuristically optimal number of rounds. In doing so, they would have taken
inefficiently high risks, since marginal losses exceeded marginal gains in later
rounds of each junctions. For instance, instead of aiming to add five or ten more
units of information value when already sitting at 70, subjects should stop and
continue with the next traffic junction to obtain twenty-five units of
information value with certainty at the same crash risk.

\subsection{Overconfidence and underconfidence}
\label{sec:overc-underc}

To determine whether a subject acts overconfidently, optimally or even underconfidently
at a certain traffic junction, the number of rounds flown is compared to the heuristics discussed before, see \eqref{eq:drift}. If the subject flew beyond the heuristically optimal number of rounds, the decision at this junction is noted as overconfident. If the subject decides to fly less rounds than optimal the junction's decision indicates underconfidence. In case the heuristic is met the subject acts optimizing.


As mentioned before, flying until the fourth node is reached
(associated with an information value of 70) is a suitable heuristic
for  a risk neutral
decision maker. This induces different strategies for the treatments. In the
closed loop system, it is optimal for the subjects to fly as many rounds as
necessary to reach the fourth node. The actual number of rounds needed depends on
the sequence of increases and non-increases   experienced, as illustrated in Fig.~\ref{fig:qualitygains}. In the open loop
system, meanwhile, flying five rounds is a suitable heuristic. As subjects have
to decide on the number of rounds to fly for all 
traffic junctions in advance, this strategy applies to all of them, as explained
before.

\begin{table}[tp]
\caption{Average shares of overconfident, underconfident and optimal decisions. Standard deviations are presented in parenthesis.}
\label{tab:average}
\begin{tabular}{| c || c | c | }
\hline
\hline
   &  Closed Loop   & Open Loop \\
\hline 
Overconfidence &  0.4410  &   0.4167 \\
degree & (.3628)  & (.3497)\\
\hline
Underconfidence &  0.3091 & 0.3833\\
degree & (.3099)	& (.3309)    \\ 
\hline
Optimizing  & 0.2521 &  0.200 \\
degree & (.2581) & (.2388)\\
\hline
\hline
\end{tabular}
\end{table}
 
\par  An indicator for overconfidence on the subject-level is created by calculating the quotient of the number of overconfident junctions and the total number of junctions per subject, since the number of traffic junctions played differed between subjects due to some drones crashing prematurely.
This relative frequency of overconfident junctions by a subject will be labeled ``overconfidence degree''. Degrees of underconfidence and optimizing behavior are
computed analogously, so all three degrees add up to 1. The average degrees of each of these behavioral tendencies are displayed in Table \ref{tab:average}.

\par
\begin{figure}[tp]
\centering
\includegraphics[scale=.35]{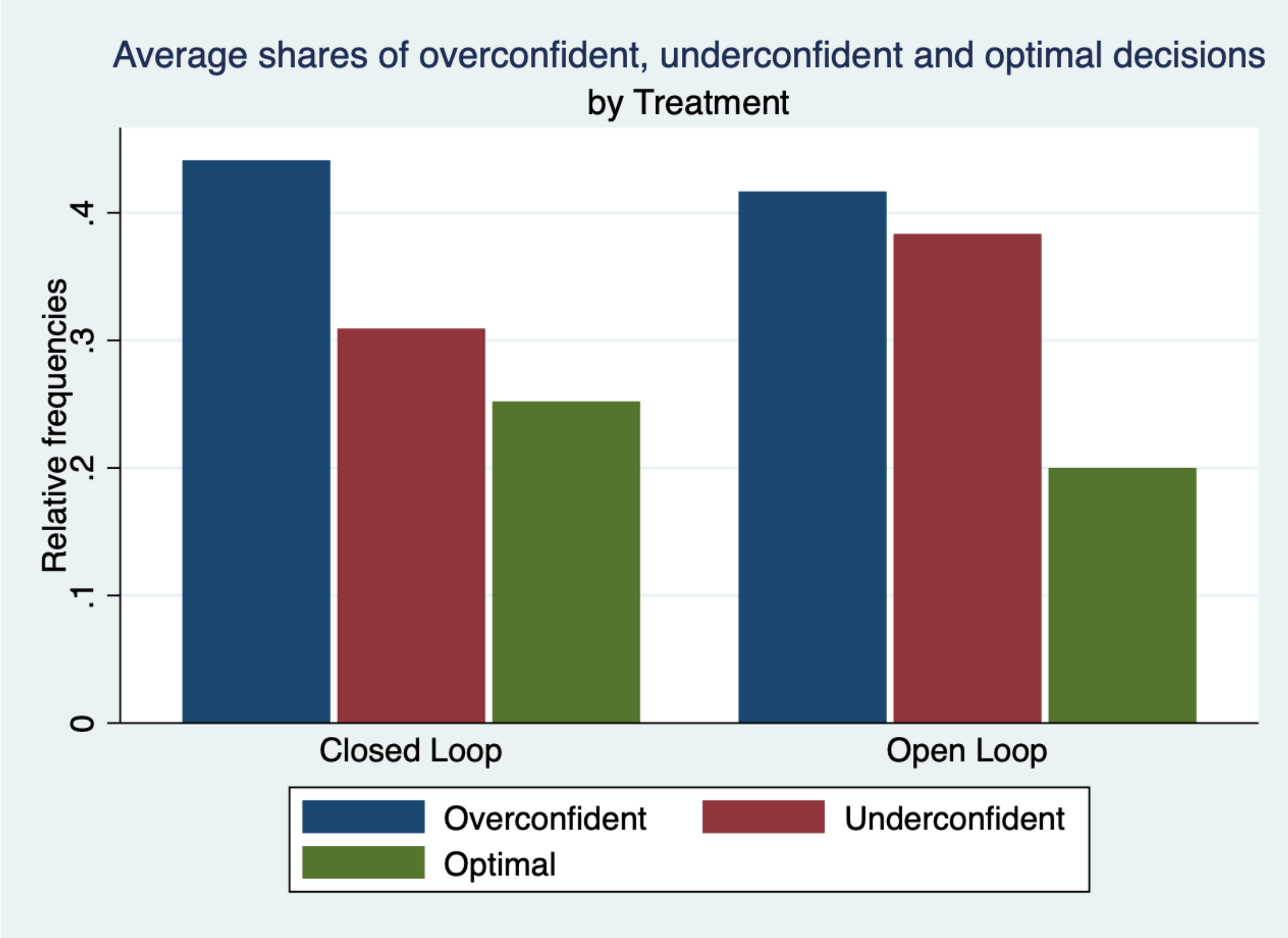} 
\caption{Average shares of overconfident, underconfident and optimal decisions}
\label{fig:averageshares}
\end{figure}

 The average overconfidence degree of 0.44 in the closed loop system means that
subjects, on average, decided overconfidently for 44$\%$ of the traffic
junctions. Comparatively, overconfidence degrees
do barely differ between the two treatments with subjects deciding overconfidently for over 42$\%$ of the junctions. The  average underconfidence degree was
higher in the Open Loop treatment. Degrees
of overconfidence and underconfidence are relatively close to each other for the
open loop, while overconfidence characterizes the predominating behavioral tendency in the
closed loop system. The optimization degree drops off in both treatments, with
subjects only deciding optimally at one fourth and one fifth of the junctions respectively. This gap also becomes apparent in the graphic
illustration of the average degrees of overconfidence, underconfidence and optimization
displayed in Fig.\
\ref{fig:averageshares}. Degrees of overconfidence did not differ significantly
by gender in either treatment.

Seven subjects in the Closed Loop and eight subjects in the Open Loop group
behaved overconfidently in every single phase they flew in, resulting in the
maximum overconfidence degree of 1. On the other side, only three subjects in
the closed loop system and two in the open loop system decided optimally for
every junction, while three subjects in each treatment acted underconfidently at every junction.  
Among those junctions in which subjects act overconfidently, they flew, on average,
2.88 (SD = 1.36) more rounds after reaching the optimal node in the
Closed Loop treatment. In the Open Loop group, subjects flew on average 1.88 (SD = 0.85) rounds beyond the optimum respectively. It appears
that those subjects, who decide to fly more rounds than optimal, they do it
clearly, i.e., usually by more than one round. Also, subjects exceed the optimal number of rounds by, on average, one more round in the closed loop  compared to the open loop system once they fly beyond the optimum. The difference is statistically significant in a Mann-Whitney U-Test ($Z=7.393, p = 0.000$) and match the observation of subjects in the Closed Loop treatment flying more rounds overall and accruing more information value in total. 

\par To gain a more precise impression about the distribution of subjects’ deviations
from optimal behavior within each treatment, we classify all subjects into
broad categories.  
Subjects that act in accordance with the   strategy of flying until reaching the
fourth node at all traffic junctions are classified as \textit{optimizing}. A subject whose 
overconfidence degree exceeds one third up to 0.5 is classified as
\textit{rather overconfident}. Once a subject flies beyond the optimal number of
rounds in more than half of the junctions (i.e., displaying an overconfidence degree
over 0.5), this subject is considered \textit{strongly
  overconfident}. Analogously, the categories \textit{rather underconfident} and
\textit{strongly underconfident} are defined using the same thresholds in terms of underconfidence degrees as the ones for overconfidence. The sixth category,
\textit{mixed}, defines scenarios, in which subjects show overconfident or
underconfident behavior in at least one junction but in less than one third of all junctions
they played  (equaling overconfidence or underconfidence degrees between 0 and
0.333). Since these subjects act optimizing in some junctions but not in all, their
behavior still contradicts the notion of a strictly optimizing agent, although no
clear dominant decision-pattern can be assigned to this behavior.  
The distribution of subjects into these categories by treatment is displayed in
Table \ref{tab:cat}.\footnote{Two subjects from the Closed Loop treatment were
  excluded, because they crashed in the very first round at the first traffic
  junction, wherefore no statement on their decision making behavior can be
  made.}

\begin{table}[tp]
\caption{Categorization of decision behavior by treatment. Relative shares referring to each treatment are presented in parenthesis.}
\label{tab:cat}
\begin{tabular}{| c || c | c || c | }
\hline
\hline
   &  Closed Loop    & Open Loop   & Total \\
\hline 
Optimal &  3  &   2 & 5 \\
               & (5.36$\%$)  & (4.26$\%$) & (4.85$\%$)\\
\hline
Rather &   5   &  6 & 10 \\
overconfident   & (8.93$\%$)  & (12.77$\%$) & (10.68$\%$)\\
\hline
Strongly &   24  &  16 & 40 \\
overconfident   & (42.86$\%$)  & (34.04$\%$) & (38.83$\%$)\\
\hline
Rather &   8 &   5 & 13 \\
underconfident  & (14.29$\%$)  & (10.64$\%$) & (12.62$\%$)\\
\hline
Strongly &   13   &  16 & 29 \\
underconfident  & (23.21$\%$)  & (34.04$\%$) & (28.16$\%$)\\
\hline
Mixed &   3   &  2 & 5 \\
	   & (5.36$\%$)  & (4.26$\%$) & (4.85$\%$)\\
	   \hline\hline
Total &   56  &  47 & 103 \\
	   & (100$\%$) & (100$\%$) & (100$\%$)\\
\hline
\hline
\end{tabular}
\end{table}

According to Table \ref{tab:cat}, only few subjects in each treatment act optimally throughout. Rather, individual behavior deviates from the optimal number of
rounds in both directions. Approximately half of the
subjects in both groups (Closed Loop: 51.79$\%$, Open Loop: 46.81$\%$) can be
considered rather overconfident or strongly overconfident. Additionally, large
shares of rather underconfident or strongly underconfident individuals (Closed Loop: 37.50$\%$,
Open Loop: 44.68$\%$) are observed as well. From the assumptions of MDPs
 and standard economic theory follow notions of individuals as
rationally optimizing decision makers that leave no room for overconfident or underconfident
behavior. Consequently, our empirical results do not  meet these theoretical
propositions. In fact, according to a Binomial Test of the data for the Closed Loop
treatment the probability of observing 29 overconfident (rather or strongly
overconfident) subjects out of a sample of 56 is virtually zero ($p < 0.000001$), given the assumption of optimizing behavior. The same result is obtained for overconfidence in the open loop 
systems. Consequently, we can fully support \textbf{H1}, indicating a strong
general tendency of individuals to act overconfidently in CPHSs, regardless of
whether receiving immediate outcome feedback or not.  
While not part of the hypothesis, we also observe the equivalent results for underconfident behavior in both treatments.

\subsection{Hot Hand Fallacy}
\label{sec:hot-hand-fallacy}
In contrast to an open loop, only a closed loop system provides the threat of sequences of immediate outcome feedback luring subjects into a hot hand
fallacy. In our experiment, a subject in the Closed Loop treatment is defined to
have a ``hot hand'' once it reaches the optimal node through obtaining
three increases in information value in a row even though the first increase was certain. This corresponds to the common perception of three repeated outcomes as a streak \cite{Carlson/Shu:2007}. In this case, an optimizing subject
would  stop flying when following the one-stage optimal strategy, since three increases are necessary to reach the optimum and
flying beyond the optimum yields a negative marginal gain. The outcome of prior
rounds should have no relevance for the subject's decision, since the rounds’
outcomes are independent from each other. If a subject decides to continue
flying though, it violates the heuristically optimal decision rule, falling
victim to the hot hand fallacy. Overall, at 84 junctions the respective subject experienced three continuous increases in
information. In 70 out of these 84 cases (83.33$\%$),
 subjects decided to fly at least one more round thereby falling victim to the hot hand fallacy. This share of hot hand fallacies among situations that pose the threat of falling victim to it is highly statistically significant: Under the assumption of subjects being optimizing Markov decision makers, the probability of such a high proportion disobeying the optimal strategy is practically zero (Binomial Test: Pr(all subjects act optimally $\mid$ 83.33$\%$ of the hot hand situations result in the hot hand fallacy): $p < 0.0000001$).
 Comparing non-overconfident and overconfident decision making between junctions with and without a
hot hand in Table
\ref{tab:hot}, we observe those subjects who experience a hot
hand situation (70 out of 84, 83.33$\%$) to have a significantly higher chance to act overconfidently at
that junction, compared to those subjects who do not reach the optimum through
three consecutive gains in information value (123 out of 260, 47.31$\%$). Running a Chi-Squared Test to test the hypothesis statistically, we obtain a highly  significant relationship between hot hand situations and subsequent overconfident behavior, i.e., the hot hand fallacy, with an error probability $p$ of virtually zero ($\chi^{2}
(df=1): 33.46$, $p = 0.000$). Consequently, \textbf{H2} is strongly
supported.

\begin{table}[tp]
\caption{Hot hand fallacy - Chi-squared table. The test captures every instance in which a subject started flying at any traffic junction in the Closed Loop treatment.}
\label{tab:hot}
\begin{tabular}{| c || c | c || c | }
\hline
\hline
   &  Not overconfident   & Overconfident  & Total \\
\hline 
\hline
No hot hand &  137  &  123  & 260 \\
\hline
Hot hand &   14 &  70  & 84 \\ 
\hline
\hline
Total &   151 & 193 & 344 \\ 
\hline
\hline
\end{tabular}
\end{table}

\begin{table}[tp]
\caption{Distribution of risk categories by treatment. Percentage values in parenthesis display relative frequencies of risk preferences in the respective treatment.}
\label{tab:risk_treat}
\begin{tabular}{| c || c | c || c | }
\hline
\hline
    &  Closed Loop  & Open Loop   &  Total\\
\hline
Risk averse & 36 & 33 &  69\\
& (63.16$\%$) & (68.75$\%$) & (65.71$\%$)\\
\hline
Risk neutral & 6 & 5 &  11\\
& (10.53$\%$) & (10.42$\%$) & (10.48$\%$)\\
\hline
Risk seeking & 6  & 4 &  10\\
& (10.53$\%$)  & (8.33$\%$) & (9.52$\%$)\\
\hline
Not identifiable & 9 & 6 &  15\\
& (15.79$\%$) & (12.50$\%$) & (14.29$\%$)\\
\hline
\hline
Total & 57 & 48 &  105\\
& (100$\%$) & (100$\%$) & (100$\%$)\\
\hline
\hline
\end{tabular}
\end{table}

\subsection{Risk Preferences}
\label{sec:risk-preferences}
Subjects' individual risk preferences are broadly classified to be risk
neutral, risk averse or risk seeking, see \cite{Dohmen/Falk/Huffman/Sunde:2010}
and Sidebar ``Risk Preferences''. Given that the presented strategy is based on
the assumption of a risk neutral decision maker, the observed distribution of
risk preferences in each treatment is displayed in Table
\ref{tab:risk_treat}.\footnote{Risk preferences of fifteen subjects in total
  could not be determined, since their choices switched between Option A and Option B three times or more.}
Consistent with the literature, the clear majority of around two thirds of the
subjects in both treatments displays risk aversion. Overall, the distribution of
risk preferences appears very similar between the treatments, with the relative
shares of each risk attitude varying only slightly between the groups. This
impression is supported by a two-sided Kolmogorov-Smirnov test, in which the
null hypothesis that the risk preference distribution in the two
treatments is not statistically different from each other cannot be rejected ($D=0.030$, $p > 0.1$). Also, degrees of overconfidence and underconfidence did not differ significantly in between-treatment comparisons for each risk attitude using Mann-Whitney U-tests.

Further, we did not observe risk seeking subjects to display significantly more overconfident behavior (Open Loop: $\varnothing = 0.5257$, Closed Loop: $\varnothing =0.175$) than risk neutral (Open Loop: $\varnothing = 0.4333$, Closed Loop: $\varnothing =0.700$) or risk averse (Open Loop: $\varnothing = 0.4140$, Closed Loop: $\varnothing =0.3818$) individuals in either treatment (Kruskal-Wallis equality of populations tests, Closed: $H = 0.450, p = 0.7984$, Open: $H = 5.437, p = 0.0660$)\footnote{Note that statistical tests on risk preferences have to be interpreted with caution, since the sample size of risk neutral and risk seeking individuals is extremely low.}. 
On the other hand, in the closed loop system, we find risk averse subjects to be significantly more underconfident, compared to the risk neutral and risk seeking subjects ($\varnothing$ underconfidence degree: $0.3480 (SD = 0.3018)$, Kruskal-Wallis equality of populations test: $H = 8.891, p = 0.0117$), while no effect was found in the open loop system.

Beside overconfidence, risk seeking preferences could in theory present an
alternative explanation for drone pilots flying beyond the optimal number of
rounds, with subjects primarily aiming to achieve a large information value
while just hoping for their drones not to crash. However, we do not find substantial differences in overconfidence conditional on risk attitude between treatments and overall mainly observe
individuals who can be classified as risk averse. This would theoretically predict
subjects to decide rather conservatively through settling for a lower number of
rounds in order to remain in the game and to protect the drone from crashing,
even at the cost of missing out higher payments from information value
gains. Since the experimental results display a large share of subjects flying
beyond the optimum, the case for individual overconfidence as the
dominant explanation for inefficient drone piloting beyond the optimum is
strengthened. Consequently, we conclude risk preferences to not being able to explain overconfident behavior although it appears to explain
the observed degrees of underconfidence).

\section{Discussion and conclusion}
\label{sec:conclusions}

Our experimental drone framework allows to observe how individuals behave when
faced with the task of piloting an UAV under risk and uncertainty, paralleling a
real-world decision-problem. Even though a closed feedback loop was identified to be
the more successful system compared to  open loop operation, we still observed
inefficient drone piloting from the vast majority of subjects. Individuals
expressed both overconfident and underconfident behavioral tendencies, regardless of receiving immediate outcome feedback. Specifically, our results indicate that immediate outcome feedback, that is originally intended to support optimal decision making, can turn out to be rather counterproductive in this regard. Overconfident decisions and consequently inefficient drone piloting can be facilitated by the hot hand fallacy as a misinterpretation of random sequences of immediate feedback on positive outcomes, since subjects fail to realize such sequences to be caused by chance and therefore being history independent. In fact, a handful of subjects in the closed loop treatment stated in the questionnaire that their decision strategy was to fly as long as they achieved steady increases in information value, trying to exploit an apparent hot hand. The fact that the possibility for this fallacy is only given in a closed loop system, presents an obvious weakness that should be considered in designing such feedback policies.

In general, the current  work exposes the human as an under observed source of errors in human-in-the-loop control systems. We thus advise researchers and practitioners to carefully account for the behavioral component in the control of cyber-physical systems and the potential problems that arise from it, besides mathematical model optimization only. In particular, our findings illustrate the impact of behavioral biases regarding effects of immediate feedback and the (miss-)understanding of history independence in chance processes.

While more information is commonly regarded to result in better decisions in
cyber-physical systems, human susceptibility for perceptual biases in response
to high information supply has to be taken into account. Therefore,
identifying an optimal quantity and frequency of feedback remains a goal for
future research. We expect a  carefully crafted intermittent feedback to be better suited
for this purpose and stress the need for an intelligent feedback design that
adapts to an individual's rationality in order provide suitable amount of
information.

Generally, considering the effect of humans in control loops more seriously
presents an important issue for research and practice. Overall, humans were shown to mostly not act optimizing in the given decision-problem in our experimental
framework, which strongly puts the Markov decision maker as an adequate
characterization of human decision making in question. Models of human decision
processes should be revisited to account for limits of cognitive capacities and
behavioral biases that result from them, in order to not jeopardize
technological accomplishments through erroneous human decisions. Otherwise,
individuals in human-in-the-loop control might take unnecessarily high risk and
render thoughtfully designed policies inefficient, as seen for highly frequent
feedback in the case of the hot hand fallacy. 
 
Lastly, our study further provides a methodological contribution to research on
CPHSs, making a first approach to incorporate insights from behavioral economics
into control engineering. Further it introduces incentivized economic
experiments as a viable option to reveal how individuals actually behave, in
contrast to how they are theoretically prescribed to behave. We present an
experimental UAV framework featuring a sequential decision-problem, with a focus
on behavioral biases in relation to feedback policies. Future experimental
research in this area may intensify efforts of incorporating various other
behavioral phenomena and stylized facts into control engineering by building
upon this framework, in order to design or test behavioral interventions that are able to
proactively counteract them.

\bibliographystyle{ieeetr}
\bibliography{litCSM}


\setcounter{equation}{0}
\renewcommand{\theequation}{S\arabic{equation}}
\setcounter{table}{0}
\renewcommand{\thetable}{S\arabic{table}}
\setcounter{figure}{0}
\renewcommand{\thefigure}{S\arabic{figure}}

\section*{Sidebars}
\subsection[Behavioral Economic Engineering]{Sidebar: Behavioral Economic Engineering}
\label{side:behavioral}
The deviation between theoretically prescribed behavior and actual human decision making is often labeled “behavioral messiness” \cite{Bolton/Ockenfels:2012} and poses a huge challenge for organizational policy makers. Experiments thereby represent a tool for bridging economic theory and real-world institutional design in order to generate practical value from applied economic science \cite{Bolton/Ockenfels:2012}. The laboratory basically acts as a “wind tunnel” for practice in order to test behavioral reactions to institutional interventions in a controlled environment, before implementing them into the real world. Projects following this approach are summed under the label of “behavioral economic engineering” \cite{Bolton/Ockenfels:2012}. Prominent examples include the implementation of optimized retirement savings plans \cite{Thaler/Benartzi:2004} and the redesign of matching algorithms for American physicians \cite{Roth/Peranson:1999}.

\subsection[Experimental research method and induced value theory]{Sidebar:
  Experimental research method and Induced Value Theory} 
\label{side:experimental}

Behavioral economic findings are often based on experimental research,
especially controlled laboratory experiments, as ``experimental control is
exceptionally helpful for distinguishing behavioral explanations from standard
ones'' \cite{Camerer/Loewenstein:2003}. Laboratory experiments, which had been
considered unfeasible for economic disciplines and privilege of natural sciences
up until the late 20th century \cite{Samuelson/Nordhaus:1985}, aim to parallel
real-world situations in laboratory settings, while abstracting from
environmental factors. This means, experiments isolate certain variables of
interest from more complex real-world contexts, while simultaneously controlling
for conditions of the subjects' economic and social environments \cite{Roth:1988}. Participants, commonly referred to as ``subjects'', work on computerized or
analogue tasks, like solving math problems \cite{Mazar/Amir/Ariely:2008} or
assembling Lego figures \cite{Ariely/Kamenica/Prelec:2008}, accompanied by
anonymously making decisions, which are observed by the experimenter. Economic
theories or concepts that the experimenter wants to test usually underlie the
decision to make in the experiment often in relation to a given task. 
A constituting factor of experiments consists in them being incentivized,
i.e., the participants are paid for solving tasks and making decisions by the
experimenter, with the concrete amount of payments depending on the respective
experimental design \cite{Roth:1988}. Subjects know about the payments they are able
to receive beforehand, as they receive instructions containing the experimental
procedure and the payoff function in the beginning of every experiment. Herein
lies the biggest advantage of experimental economic research: actual human
behavior can be observed with the subjects’ actions having actual monetary
consequences for their payoffs in the end, so they have to ``put their money
where their mouth is''. This is not the case for other research methods used in
economic or social science like surveys or scenario studies, in which
participants only state how they would behave in certain situations, while no
information is gained whether they would actually behave the way they stated
when faced with the decision-problem in reality. The problems of
intention-behavior gaps \cite{Carrington/Neville/Whitwell:2010} and giving socially
desirable answers \cite{Nederhof:1985} are well-known. 
\par As a basic principle of economic experiments, there is no deception by the
experimenter. Subjects will be asked to do exactly what is stated as their task
in the instructions they receive before the experiment and they will be paid
exactly according to the payoff function given in the instructions as well
\cite{Weimann/BrosigKoch:2019}. Experiments can be used for various purposes
like testing economic theories and theoretical equilibria, testing policies and
environments, establishing phenomena, stylized facts and new theories or
deriving political recommendations \cite{Roth:1988}. Furthermore, experimental
evidence can be replicated and re-evaluated by other experimenters using the
same experimental setup and instructions, which are usually published alongside
the results \cite{Charness:2010}. 

In order to experimentally test certain theories or economic interventions,
subjects are randomly assigned into one control group and (at least) one
treatment group, in order to avoid (self-)selection biases, with a treatment
being the intervention that is going to be tested. Differences between control
group and treatment group must thereby be reduced to exactly one factor, namely
the so-called treatment variable, i.e., the intervention that the experimenter wants to
evaluate \cite{Weimann/BrosigKoch:2019}. 
Herein lies another major advantage of controlled economic experiments as they
allow causal inferences. Due to the experimental environment being held
constant and environmental factors being controlled for all groups, they cancel
out through comparative statics between groups. Therefore, outcome in a
treatment group compared to a control group can only be caused by an exogenous
change, i.e., by the treatment intervention itself since it represents the only
factor that differs between groups \cite{Pearl:2003}, \cite{Camerer/Loewenstein:2003}, \cite{Weimann/BrosigKoch:2019}. 
 
\par As mentioned above, a constituent factor of economic experiments consists
in the provision of decision-dependent monetary incentives. Economic experiments
underlie the general assumption that any kind of utility an individual
experiences can be expressed by an equivalent monetary incentive. Through these
incentives, the experimenter is able to induce a certain utility function $U(x)$
onto the subjects, in order to neutralize the subjects’ inherent preferences for
the duration of the experiment. This ``induced value theory'' was originally 
introduced by Nobel Memorial Prize in Economic Sciences laureate Vernon Smith \cite{Smith:1976}. 
According to him, subjects should derive all utility from the monetary incentives provided in the experiment
itself, during their participation.

\par While the internal validity of experiments is usually recognized to be very
strong, the external validity of these experiments, i.e., the transferability of
their results to the world outside the laboratory, the external validity, is
commonly discussed and criticized by opponents of this methodology
\cite{Slovic/Fischhoff/Lichtenstein:1978}, \cite{Weimann/BrosigKoch:2019}.   
For instance, in the context of perceptional and judgmental biases, laboratory
settings are accused of changing the environment in which a human makes
efficient decisions to an artificial one. However, it can be argued the other
way around in that laboratory experiments show people at their best through
providing all information necessary and eliminating distractions
\cite{Slovic/Fischhoff/Lichtenstein:1978}. If people fall to biases in this
isolated, save environment, they will likely do so outside of it as well.

\subsection[Risk Preferences ]{Sidebar: Risk Preferences }
\label{side:riskpreferences}

\begin{figure}[tp]
\centering
\includegraphics[scale=0.3]{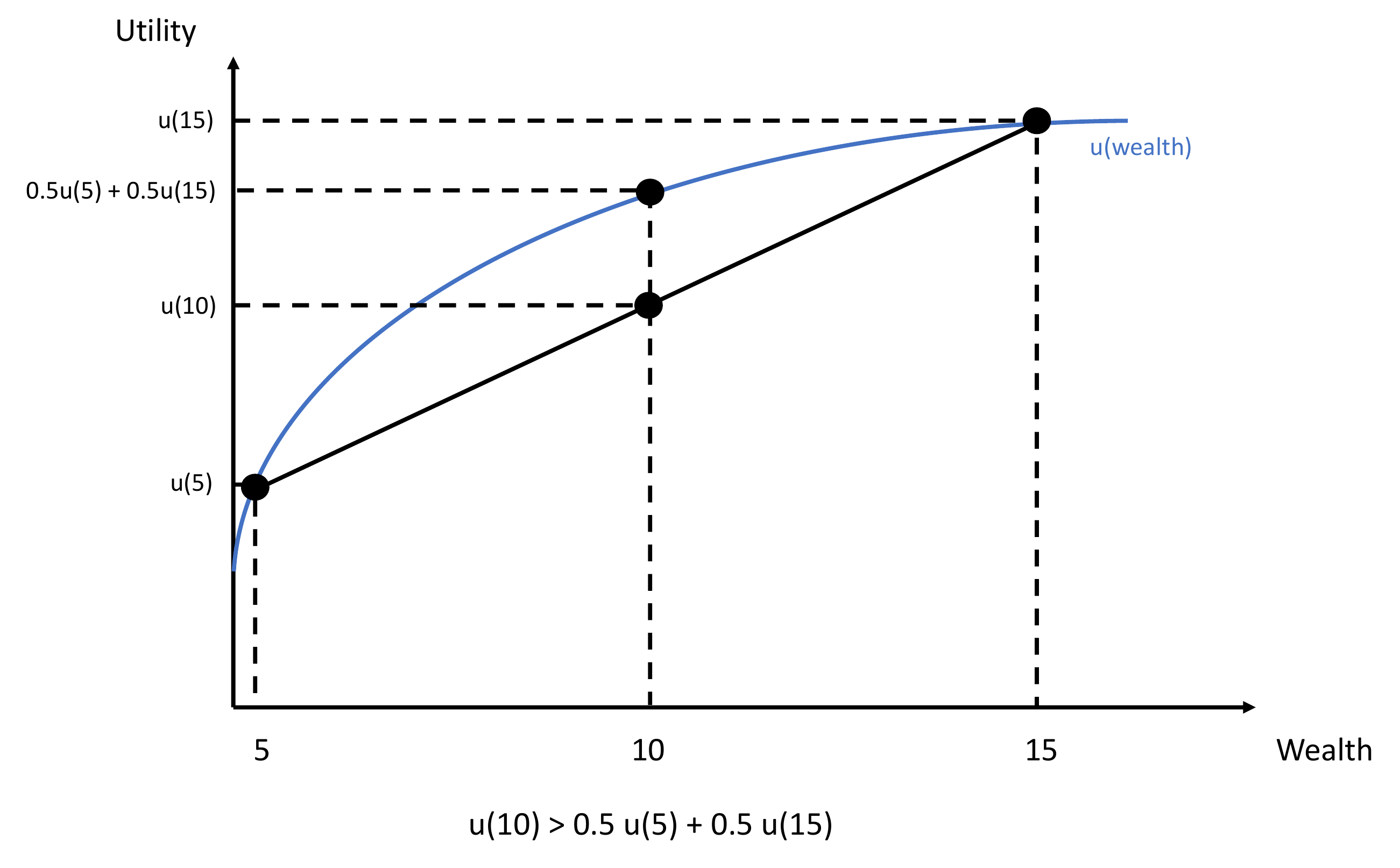} 
\caption{Example for a risk averse individual's utility function, adapted from \cite{Varian:2006} }
\label{fig:riskaverse}
\end{figure}

Individual risk preferences represent a factor commonly influencing human
behavior that should necessarily be considered in order to explain individual
decision making under risk and uncertainty. Human risk preferences are generally
differentiated into risk averse, risk neutral and risk seeking behavior. They
are usually determined by the individual's subjective evaluation of a
probabilistic payment from a lottery compared to a certainty equivalent, i.e., a safe fix payment \cite{Holt/Laury:2002}, \cite{Dohmen/Falk/Huffman/Sunde:2010}. 
If an individual subjectively judges the utility from the expected value $E$ of
a payment $X$ as higher than the expected utility of such payment, namely
$u(X)$, then this individual is classified as risk averse  \cite{Pratt:1964}: 
$$u(E[X])>E[u(X)]$$
As a consequence of Jensen's inequality, a risk averse individual's personal
utility function is concavely shaped. This means that the individual is willing to
sacrifice the chance of higher but uncertain payments in order to receive a
smaller but certain payment, the certainty equivalent, thereby paying the
so-called insurance premium  \cite{Pratt:1964}. In an example by
\cite{Varian:2006} displayed in Fig.\ \ref{fig:riskaverse}, an individual with
\euro{10} at his disposal has to decide whether or not to participate in a
lottery in which \euro{5} are gained with a probability of 50$\%$ and \euro{5}
are lost with a probability of 50$\%$. For a risk averse individual the utility
of the expected value $u(10)$ is greater than the expected utility of wealth
$0.5u(15) + 0.5u(5)$. 
\par
\begin{figure}[tp]
\centering
\includegraphics[scale=0.3]{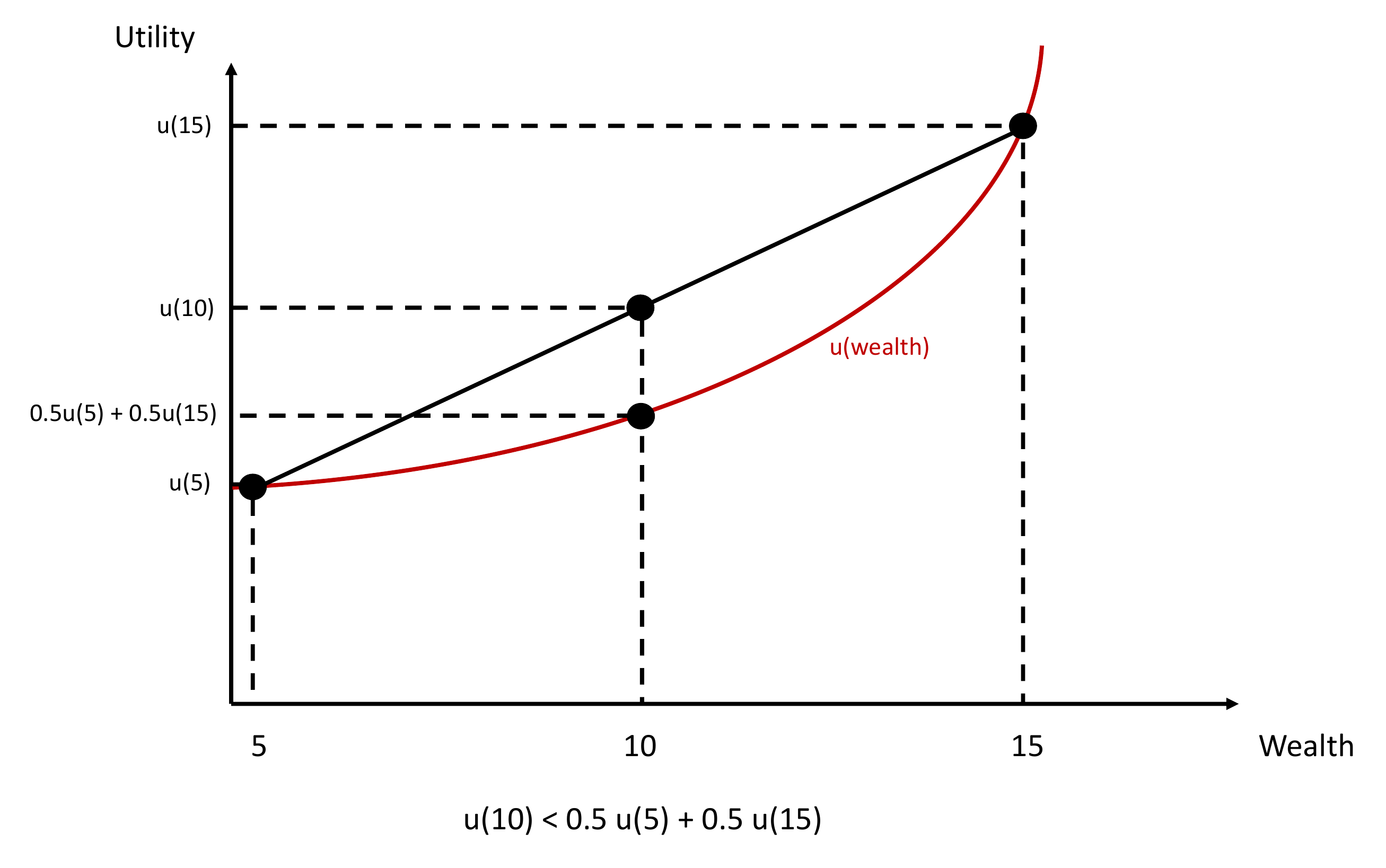} 
\caption{Example for a risk seeking individual's utility function, adapted from \cite{Varian:2006} }
\label{fig:riskseeking}
\end{figure}

 The other way around, if an individual derives less utility from the
expectated value of $X$ than the $X$’s expected utility $$u(E[X])<E[u(X)],$$
then  this individual is considered risk seeking. Its subjective utility
function is convex (see the analogous example in Fig.\
\ref{fig:riskseeking}). Risk seeking individuals are willing to give up the
certainty equivalent in order to get the chance to play the lottery through
which they may obtain a higher payment, i.e., they favour gambling. The payment
obtainable beyond the certainty equivalent through playing the lottery is called
risk premium. In case an individual is indifferent between both alternatives,
the individual can be considered risk neutral, with risk and insurance premium
cancelling out \cite{Pratt:1964}.  

Standard economic theories, in the tradition of the \textit{homo oeconomicus}, are
usually based upon the explicit or implicit assumption of individual risk
neutrality, reflecting notion of individuals as rationally optimizing (``Markovian'') decision makers  \cite{Jaquette:1976}. Contrarily, behavioral economic research usually observes risk averse decision makers, in the laboratory  \cite{Harrison/List/Towe:2007}, \cite{Eckel/Grossman:2008}, as well as in the field  \cite{Harrison/Lau/Rutstrom:2007}.
 
\par The degree of human risk aversion was found to increase with the size of incentives being at stake \cite{Kachelmeier/Shehata:1992}, \cite{Holt/Laury:2002}, contradicting theories of constant relative risk aversion \cite{Holt/Laury:2002}. However, individuals are observed to act clearly risk averse even only with low incentives at stake, wherefore models that assume risk neutral human behavior do not adequately represent reality and may lead to biased inference \cite{Holt/Laury:2002}. Regardless of this quite clear evidence, only sporadic effort exists in order to incorporate risk preferences other than risk neutrality into Markov Decision Processes  \cite{Howard/Matheson:1972}, \cite{Majumdar/Singh/Mandlekar/Pavone:2017}.

\subsection[Mulitiple Price List Format ]{Sidebar: Mulitiple Price List Format }
\label{side:multiple}

Multiple price lists formats offer a simple but effective tool to reveal
subjects' actual risk preferences. The basic concept was designed with the goal of testing the standard economic assumption of risk
neutrality as well as the behavioral assumption of constant relative risk
aversion. The list consisted of a table with two
columns that represented two lotteries the subjects had to choose between for
ten rows featuring varying probabilities of winning certain amounts of money
that remained constant for all rows \cite{Holt/Laury:2002}. As the original
price list contained monetary values with two decimals and
varying probabilities that had to be weighed, it can be assumed that it was
difficult for subjects to compare the different options. The subsequent work \cite{Dohmen/Falk/Huffman/Sunde:2010} designed a
more clear-cut version of the 
multiple price list (see Table \ref{tab:multiple}) in which subjects had to
choose between a safe payoff that increased by one Euro with each row and a fair
lottery that remained constant for all rows.  

\par Subjects are asked to indicate for each row whether they prefer the offered
safe payment (Option A) or the offered fair lottery (Option B). The safe payment
equals the number of the respective row minus one (in Euro), starting at \euro{0}
in the first row and increasing up to \euro{19} in the twentieth row. The
lottery is the same in each row, containing a $50\%$ chance to win \euro{30} and
a $50\%$ chance of winning nothing. A subject's risk attitude can be determined
by observing the row in which his choice changes from B to A. When continuously
choosing B until row 16 and continuously choosing A afterwards, a subject is
classified as risk neutral, as both options yield the exact same expected
payoffs of \euro{15} in row 16. At this point a risk neutral subject is
indifferent between the receiving the safe payment and playing the lottery. If
the change from B to A occurs before the sixteenth row, then a subject is classified
as risk averse, as it is willing to sacrificing the chance to play a lottery
that yields a higher expected payoff in favor of a safe payment of a lower
amount than such expected payoff. If the change occurs after the sixteenth row,
a subject is regarded as risk seeking. Contrarily to risk averse subjects,
risk seeking subjects are willing to sacrifice a safe payment in order to be
able to play a lottery that yields an expected payoff that is lower than safe
payment but yields a higher maximum payoff. 

\par After  making their choices for every row, a fraction of the
participant group (e.g., every seventh \cite{Dohmen/Falk/Huffman/Sunde:2010} or fifteenth
\cite{Djawadi/Fahr:2013} subject) is chosen at random and actually paid for one
randomly chosen row of the table, according to their respective choice made in
this row. This procedure incentivizes subjects to choose in accordance with
their true preferences in each row, since they would forfeit expected utility if
they did not, thereby revealing them to the experimenter \cite{Dohmen/Falk/Huffman/Sunde:2010}. Therefore, the multiple price list format presents an incentive-compatible instrument to measure individual risk preferences.

\begin{table}[t]
\caption{Multiple price list by \cite{Dohmen/Falk/Huffman/Sunde:2010} }
\label{tab:multiple}
\begin{tabular}{| c | l | l | }
\hline
\hline
& \textbf{Option A} & \textbf{Option B}\\
\hline
1) & \euro{0}  safe  & \euro{30} with a prob.\ of $50\%$, \euro{0} with a prob.\ of $50\%$ \\
	\hline
2) & \euro{1}  safe  & \euro{30} with a prob.\ of $50\%$, \euro{0} with a prob.\ of $50\%$  \\
\hline
3) & \euro{2}  safe  & \euro{30} with a prob.\ of $50\%$, \euro{0} with a prob.\ of $50\%$  \\
\hline
4) & \euro{3}  safe  & \euro{30} with a prob.\ of $50\%$, \euro{0} with a prob.\ of $50\%$  \\
\hline
\vdots &\vdots &\vdots  \\
\hline
19) & \euro{18}  safe  & \euro{30} with a prob.\ of $50\%$, \euro{0} with a prob.\ of $50\%$  \\
\hline
20) & \euro{19}  safe  & \euro{30} with a prob.\ of $50\%$, \euro{0} with a prob.\ of $50\%$  \\
\hline
\hline
\end{tabular}
\end{table}


\section{Authors' Biographies}

\textbf{\textit{Marius Protte}} is a
Ph.D. student at the Chair of Corporate Governance at Paderborn University's Faculty of Business Administration and Economics. He
received his M.Sc.\ degree in Business Administration in 2019 from Paderborn
University. His research interests are in Behavioral Economics and Business
Ethics and Human-Machine Interaction.

\textbf{\textit{Ren\'e Fahr}} received his diploma
(M.Sc.) degree in Economics in 1998 from the University of Bonn. In
2003 he was 
awarded the Ph.D.\ degree jointly from the 
 University of Bonn and the London School of Economics as part of the European Doctoral Program in quantitative economics.
 Since 2008 he is the head of Chair of Corporate Governance at Paderborn University and Research Director of the BaER-Lab at the Faculty of
 Economics and Business Administration. Since 2019 he is Vice-President for Technology Transfer and member of the University's executive board
 and since 2020 he is head of the research group Behavioral Economic Engineering and Responsible Management at the Heinz Nixdorf Institute (Paderborn University).
 His doctoral studies were supported by full scholarships by the German Science Foundation and German Academic Exchange Service.
 Since 1998 he was affiliated with the IZA in Bonn as a resident research affiliate. 
 From 2003 till 2008 he was a post-doctoral researcher in personnel economics at the University of Cologne.
  His research interests are in Behavioral Economics, Behavioral Business Ethics
  and Human-Machine Interaction. 

\textbf{\textit{Daniel Quevedo}} is professor of cyberphysical systems in the School of Electrical Engineering and Robotics, Queensland University of Technology (QUT), Brisbane, Australia.  
Before joining QUT, he established and led  the Chair in Automatic
Control 
at Paderborn University, Germany. In 2000, he received Ingeniero Civil Electr\'onico
and M.Sc.\ degrees from   Universidad
T\'ecnica Federico Santa Mar\'{\i}a, Valpara\'{\i}so, Chile,  and in 2005   the Ph.D.\ degree from the University of Newcastle in Australia.
He  received the IEEE Conference on Decision and
Control Best Student Paper Award in 2003 and was also a finalist  in 
2002.  In
2009 he was awarded a five-year  Research Fellowship from the Australian
Research Council. He is co-recipient of the  2018 IEEE Transactions on Automatic Control George S.~Axelby Outstanding Paper Award. 
Prof.~Quevedo currently serves as Associate Editor for the \textit{IEEE Control Systems Magazine} and in  the Editorial Board of the \textit{International Journal of Robust and Nonlinear Control}.  From 2015--2018 he was Chair of the IEEE Control Systems Society \textit{Technical Committee on 
  Networks \& Communication Systems}. 
 His  research
interests are in  control of networked systems and of power
converters.  
 
\end{document}